\documentclass{article}
\usepackage{maa-monthly}
\usepackage{url}

\theoremstyle{theorem}
\newtheorem{theorem}{Theorem}

\theoremstyle{definition}

\begin{document}

\title{An identity for the infinite sum $$\sum_{n=0}^\infty\frac{1}{(n!)^3}$$}
\markright{Submission}
\author{Hassan Jolany}

\maketitle

\begin{abstract} In this short note, we give an identity for the alpha function $$\alpha(x,s)=\sum_{n=0}^\infty\frac{x^n}{(n!)^s},$$ where $s\in\mathbb N, \; x\in \mathbb R$, in the case $s=3$. 
\end{abstract}

\noindent

\section{Introduction} We know that when $s=1$, then the alpha function is the exponential function $\alpha(x,1)=e^x$. The question is that what about when $s\geq 2$ ?. We have the following nice differential equation for the alpha function \cite{Eguether}. In fact, it is known that the alpha function is the solution of the following differential equation.

$$\sum_{k=1}^s\sigma_s^kx^{k-1}y^{(k)}-y=0$$

where $\sigma_s^k$ is the Stirling numbers of the second kind which satisfies in the following generating function,
$$\frac{1}{k!}(e^x-1)^k=\sum_{s=k}^\infty\sigma_s^k\frac{x^n}{n!}$$

Hence by using the language of generalized hypergeometric function, we get \cite{Alpha}, $$\sum_{n=0}^\infty\frac{1}{(n!)^3}=_0 F_2(;1, 1;1)=\sim 1.1297
$$

To deal with this question, we first give the following theorem which is the direct consequence of Bessel-Parseval identity\cite{Konrad}. 

\begin{theorem}Let we have the two following entire series, with real the coefficients, $$f(x)=\sum_{n=0}^\infty a_nx^n \;\text{and}\; g(x)=\sum_{n=0}^\infty b_nx^n  $$ with respectively non-zero radius $R$ and $R'$. Take $$h(x)=\sum_{n=0}^\infty a_nb_n x^n$$, then if $-R<u<R$ and $-R'<v<R'$, then we have

$$h(uv)=\frac{1}{2\pi}\int_0^{2\pi}f(ue^{it})g(ve^{-it})dt$$
\end{theorem}

In the Theorem 1, if we take $f(x)=g(x)=e^x$, and $u=x, v=1$,  then we get the following amazing identity,
$$\sum_{n=0}^\infty\frac{x^n}{(n!)^2}=\frac{1}{2\pi}\int_0^{2\pi}e^{(x+1)\cos t}\cos((x-1)\sin t)dt$$

Hence by taking $x=1$ we get

$$\sum_{n=0}^\infty\frac{1}{(n!)^2}=\frac{1}{2\pi}\int_0^{2\pi}e^{2\cos t}dt$$

But we know, 

\begin{equation}
 \frac{1}{2\pi} \int\limits_{-\pi}^{\pi} e^{\left(a \cos x + b \sin x\right)} dx = I_0\left(\sqrt{a^2+b^2}\right),
\end{equation}
where $I_0$ is the modified Bessel function of the first kind.

\section{The main Result}In this section we give an integral identity for the alpha function when $s=3$. We take $f(x)=e^x$,  $g(x)=\frac{1}{2\pi}\int_0^{2\pi}e^{(x+1)\cos t}\cos((x-1)\sin t)dt$, and $u=x, v=1$. Hence, from Theorem 1, we have

$$A=\sum_{n=0}^\infty\frac{x^n}{(n!)^3}=\frac{1}{4\pi^2}\int_0^{2\pi}\int_0^{2\pi}e^{xe^{is}}e^{(e^{-is}+1)\cos t}\cos \left((e^{-is}-1)\sin t\right)dsdt$$

But, we have $$e^{xe^{is}}e^{(e^{-is}+1)\cos t}=\left(e^{x\cos s+\cos s\cos t+\cos t}\right)\left(\cos (x\sin s-\sin s\cos t)+i\sin(x\sin s-\sin s\cos t)\right)$$

and also,
$$\cos \left((e^{-is}-1)\sin t\right)=\cos(\cos s\sin t-\sin t)\cosh(\sin s\sin t)+i\sin(\cos s\sin t-\sin t)\sinh(\sin s\sin t)$$

Hence

\begin{align*}
A&=\frac{1}{4\pi^2}\int_0^{2\pi}\int_0^{2\pi}\left((e^{x\cos s+\cos s\cos t+\cos t})(\cos (x\sin s-\sin s\cos t))\right)\cos(\cos s\sin t-\sin t)\cosh(\sin s\sin t)dsdt\\
&-\frac{1}{4\pi^2}\int_0^{2\pi}\int_0^{2\pi}\left(e^{x\cos s+\cos s\cos t+\cos t}\right)\sin(x\sin s-\sin s\cos t)\sin(\cos s\sin t-\sin t)\sinh(\sin s\sin t)dsdt
\end{align*}

and this gives an identity for the alpha function in the case $s=3$.


\begin{thebibliography}{9}
\bibitem{Konrad} 
Knopp, Konrad,
\textit{Theory and Application of Inﬁnite Series}. Dover Publications.(1990) 

\bibitem{Eguether} 
Gerard Eguether, 
\textit{Nombres de Stirling et Nombres de Bell}. 
\\\texttt{http://www.iecl.univ-lorraine.fr/~Gerard.Eguether/zARTICLE/BF.pdf}

\bibitem{Alpha} 
Wolfram Alpha,\\\texttt{http://www.wolframalpha.com/}

\end{thebibliography}
\end{document}